\documentclass[11pt]{amsart}
\theoremstyle{plain}
\newtheorem{Theorem}{Theorem}
\newtheorem{Lemma}{Lemma}
\begin{document}
\title{Cordial Deficiency}
\author{Adrian Riskin}
\address{Department of mathematics\\
Mary Baldwin College\\
Staunton, Virginia 24401\\
USA}
\email{ariskin@mbc.edu}
\thanks{I would like to note the contributions of my colleagues W. Michael Gentry and John Ong,
who have provided me with
an ongoing education in some of the topics treated in this paper.}
\begin{abstract}
We introduce two new measures of the noncordiality of a graph.
We then calculate the values of these measures for various families of
noncordial graphs.  We also determine exactly which of the M\"obius ladders
are cordial.
\end{abstract}

\maketitle

\baselineskip=20pt

\section{Introduction and definitions}

Cahit [1] introduced cordial graph labelings twenty years ago, and they have remained the focus of
a steady if scant stream of papers, for a summary and listing of which see
Gallian [2].  We allow graphs to have multiple edges but not loops.  A
\textit{binary labeling} of a graph $G$ is a function $f:V(G) \rightarrow \{
0,1 \}$.  Two real numbers $x$ and $y$ are \textit{roughly equal} if $0 \leq
|x-y| \leq 1$.  A binary labeling is said to be \textit{friendly} if
$|f^{-1}(0)|$ is roughly equal to $|f^{-1}(1)|$.  A binary labeling $f$ of $G$
induces a binary labeling $f_{e}:E(G) \rightarrow \{ 0,1 \}$ by
$f_{e}(uv)=f(u)+f(v)$, where the sum is calculated modulo 2.  A friendly
labeling of $G$ is \textit{cordial} when $|f_{e}^{-1}(0)|$ is roughly equal to
$|f_{e}^{-1}(1)|$.  After the tradition of those who author new methods of
graph labeling, Cahit's seminal paper contains a number of theorems concerning
the cordiality and noncordiality of various families of graphs, the most
salient of which to our purpose is:

\begin{Theorem} The complete graph $K_{n}$ is cordial if and only if $n \leq 3$.
\end{Theorem}

We introduce two measures of the noncordiality of a graph.  These were inspired by Kotzig's and Rosa's notion of
\textit{edge-magic deficiency} [3].  Note that every friendly labeling of a
graph $G$ can be made into a cordial labeling of a graph $G^{\prime}$ by
adding no more than $|f_{e}^{-1}(0)-f_{e}^{-1}(1)|-1$ edges between
appropriate pairs of vertices so that the number of edges labeled $0$ becomes
roughly equal to the number of edges labeled $1$ in the augmented graph
$G^{\prime}$.  The minimum number of edges, taken over all friendly labelings
of $G$, which it is necessary to add in order that $G^{\prime}$ become cordial
is the \textit{cordial edge deficiency} of $G$, denoted by $ced(G)$.

If it is possible to find a binary labeling of $G$ so that $|f_{e}^{-1}(0)|$
and $|f_{e}^{-1}(1)|$ are roughly equal, then $G$ can be made cordial by
adding no more than $|f^{-1}(0)-f^{-1}(1)|-1$ vertices, labeled appropriately.
 The minimum number of vertices, taken over all such binary labelings of $G$,
which it is necessary to add in order that $G^{\prime}$ become cordial is the
\textit{cordial vertex deficiency} of $G$, denoted by $cvd(G)$.  If there are
no such binary labelings of $G$ we say that $G$ is \textit{strictly
noncordial}, and write $cvd(G)=\infty$.

\section{Results}

Our first two theorems have to do with the cordial deficiencies of the complete
graph.

\begin{Theorem}
The cordial edge deficiency of $K_{n}$ for $n>1$ is $\left \lfloor \frac{n}{2}
\right \rfloor - 1$.
\end{Theorem}

\noindent \textbf{\textit{Proof:}} Let $f$ be a friendly labeling of $K_{n}$ and
suppose that $n=2j$.  Then $|f^{-1}(0)|=|f^{-1}(1)|=j$.  Thus
$|f_{e}^{-1}(0)|=2{j \choose 2}=j^{2}-j$ and $|f_{e}^{-1}(1)|=j^{2}$.  The
difference between them is $j$, and so $cev(K_{n})=j-1$.  A similar
calculation yields the result when $n=2j+1$. \hfill $\square$

\begin{Theorem}
The cordial vertex deficiency of $K_{n}$ is $j-1$ if $n=j^{2} + \delta$, where
$\delta \in \{-2, 0, 2 \}$.  Otherwise $K_{n}$ is strictly noncordial.
\end{Theorem}

\noindent \textbf{\textit{Proof:}} Let $f$ be a binary labeling of $K_{n}$. 
Suppose that $|f^{-1}(0)| = \ell$, so that $|f^{-1}(1)| = n- \ell$.  Then
$|f_{e}^{-1}(0)| = {\ell \choose 2} + {n-\ell \choose 2} $ and
$|f_{e}^{-1}(1)|=\ell (n-\ell)$.  The difference between these is $\left | 2
\ell^{2} - 2n\ell + {n \choose 2} \right |$.  This equals $0$ if and only if
$\ell = \frac{n \pm \sqrt{n}}{2}$ and it equals $1$ if and only if 
$\ell =\frac{n \pm \sqrt{n \pm 2}}{2}$.  Since $\ell$ is a whole number, we
must have $n=j^{2}+\delta$ for $\delta \in \{-2, 0, 2 \}$ in order to have
$|f_{e}^{-1}(0)|$ roughly equal to $|f_{e}^{-1}(1)|$.  If $\delta = 0$ we have
$\ell = \frac{j^{2} \pm j}{2}$, in which case the difference between
$|f^{-1}(0)|$ and $|f^{-1}(1)|$ is $j$, so that $cvd(K_{n})=j-1$.  Similar
calculations yield the result in the other two cases.  \hfill $\square$\\

\smallskip
We next consider the cordiality of M\"obius ladders.  The M\"obius ladder $M_{k}$
consists of the cycle $C_{2k}$ (the \textit{canonical $2k$-cycle}) with $k$ additional
edges (the \textit{cross-edges}) joining opposite pairs of vertices.  These graphs have a
natural grid-like embedding into the M\"obius strip, from which they take their name, and 
of which it is fruitful to think while reading the proofs which follow.

\begin{Lemma}
If $k \equiv 2\pmod 4$ then $M_{k}$ is not cordial.
\end{Lemma}

\noindent \textbf{Proof:} Suppose $k=4n+2$.  Then $M_{k}$ has $8n+4$ vertices and $12n+6$ 
edges.  If $f$ is a cordial labeling of $M_{k}$ we have 
\begin{equation*}
\sum_{uv \in E(M_{k})} f_{e}(uv) = 6n+3 \equiv 1 \pmod{2}
\end{equation*}
Furthermore, 
\begin{equation*}
\sum_{uv \in E(M_{k})} f_{e}(uv) = 3\sum_{v \in V(M_{k})} f(v) = 12n+6 \equiv 2 \pmod 2
\end{equation*}
and we have obtained a contradiction. \hfill $\square$

\begin{Theorem}\label{cordial}
If $k \not\equiv 2 \pmod 4$ and $k \geq 3$ then $ M_{k}$ is cordial.
\end{Theorem}

\noindent \textbf{Proof:} The cordiality of $M_{3}$ follows from the labeling of the vertices of the
canonical 6-cycle with 1, 1, 0, 1, 0, 0 in this order.  That of $M_{4}$ from 1, 1, 0, 1, 1, 0, 0, 0, and that
of $M_{5}$ from 1, 1, 1, 1, 0, 1, 0, 0, 0, 0.  Now, if we have a cordially labeled $M_{k}$ which has
a cross-edge $uv$ with $f(u)=f(v)=1$ then it is possible to separate the labeled graph along this edge, 
obtaining a (non-cordially) labeled copy of $P_{2} \times P_{k}$ in the process, and do the same with
a cordially labeled copy of $M_{4}$, which can then be grafted into the modified $M_{k}$ with an 
appropriate twist, yielding a cordially labeled copy of $M_{k+4}$.  Since all three of the labelings given
above have such a cross-edge, the result follows by induction. \hfill $\square$

\begin{Theorem}
If $k \equiv 2 \pmod 4$ and $k \geq 6$ we have $ced(M_{k})=cvd(M_{k})=1$.
\end{Theorem}

\noindent \textbf{Proof:} Consider the friendly labeling of $M_{6}$ obtained by assigning the labels
1, 1, 1, 1, 1, 0, 1, 0, 0, 0, 0, 0 in this order to the vertices of the canonical 12-cycle.  This labeling has
$|f_{e}^{-1}(0)|=10$ and $|f_{e}^{-1}(1)|=8$.  Hence $ced(M_{6})=1$.  The fact that $ced(M_{k})=1$
follows by induction exactly as in the proof of Theorem \ref{cordial}.  The result for $cvd(M_{k})$ follows
by the same method using the binary labeling 1, 1, 1, 0, 1, 0, 1, 1, 0, 0, 1, 0. \hfill $\square$

The \textit{wheel graph} $W_{n}$ is obtained from the cycle $C_{n}$ (the \textit{canonical n-cycle}) by
adding another vertex (the \textit{central vertex}) and joining it to the $n$ vertices of the canonical 
$n$-cycle.  The $n$ edges incident with the central vertex are called \textit{central edges} and the other $n$
edges are called \textit{cycle edges}.  Cahit [1] showed that $W_{k}$ is cordial if and only if 
$n \not \equiv 3 \pmod 4$ and that $C_{n}$ is cordial if and only if $n \not \equiv 2 \pmod 4$.

\begin{Theorem}
If $n \equiv 3 \pmod 4$ then $ced(W_{k})=cvd(W_{k})=1$.
\end{Theorem}

\noindent \textbf{Proof:} If $n=4k+3$ then $|V(W_{n})|=4k+4$ and $|E(W_{n})|=8k+6$.  Let
$f$ be a friendly labeling of $W_{n}$.  Note that we may assume without loss of generality that 
if $w$ is the central vertex then $f(w)=0$.  This leaves $2k+1$ cycle vertices labeled 0 and $2k+2$
labeled 1, which in turn yields $2k+1$ central edges labeled 0 and $2k+2$ labeled 1.  Since $C_{n}$ is
cordial, it is possible to arrange the vertex labels on the canonical $n$-cycle so that $|f_{e^{*}}^{-1}(0)|$
is roughly equal to $|f_{e^{*}}^{-1}(1)|$, where $f_{e^{*}}$ represents the restriction of $f_{e}$ to 
cycle edges.  This implies that $|f_{e^{*}}^{-1}(0)|=2k+1$ and $|f_{e^{*}}^{-1}(1)|=2k+2$, since 
otherwise we'd have a cordial labeling of $W_{n}$.  Therefore $ced(W_{n})=1$.

Now, if we begin with a cordial labeling of the canonical $n$-cycle then we may assume without loss of
generality that we have $2k+1$ cycle vertices labeled 0 and $2k+2$ labeled 1.  As above this implies 
that we have $2k+1$ cycle edges labeled 0 and $2k+2$ labeled 1.  Hence if we label the central vertex
with 1, we obtain a binary labeling $f$ in which $|f_{e}^{-1}(0)|=|f_{e}^{-1}(1)|=4k+3$, and in which
$|f^{-1}(0)|=|f^{-1}(1)|-2$.  The result follows immediately. \hfill $\square$

\section{References}
\begin{enumerate}

\item Cahit, I.; Cordial graphs: a weaker version of graceful and harmonious
graphs.  Ars Combin. 23(1987) 201-207.

\item Gallian, J. A.; A dynamic survey of graph labeling.  Electronic J.
Combin. DS6. http://www.combinatorics.org/Surveys/index.html

\item Kotzig, A. and Rosa, A.; Magic valuations of finite graphs.  Canad.
Math. Bull. 13(1970) 451-461.

\end{enumerate}

\end{document}